\documentclass[11pt]{article}

\usepackage{amsmath,amssymb,amsthm}
\usepackage{geometry}
\usepackage{enumerate}
\usepackage{hyperref}
\newtheorem{theorem}{Theorem}
\newtheorem{lemma}[theorem]{Lemma}
\newtheorem{definition}[theorem]{Definition}

\title{Logical Undefinability of the Generalized Collatz Transition Relation in B\"uchi Arithmetic}
\author{Madhav Dhiman \thanks{Primary/Corresponding Author: mmdhiman09@gmail.com} \and Rohan Pandey \thanks{rpande@uw.edu}}
\date{\today}

\begin{document}

\maketitle

\begin{abstract}
Let $q$ be an odd prime and let $d$ be an odd integer. We show that the arbitrary-step transition relation of the generalized Collatz map $T_{q,d}$ is not first-order definable in Base-2 B\"uchi Arithmetic ($BA_2$). We do this by demonstrating that if the transition relation were definable, the exponential set $P_q = \{q^y : y \in \mathbb{N}\}$ would also be definable in $BA_2$. Since $P_q$ is strictly non-semilinear, this yields a direct contradiction with the Cobham--Sem\"enov theorem. Consequently, we demonstrate that no finite automaton reading base-2 representations can recognize this transition relation.
\end{abstract}

\section{Introduction}

Let $q$ be an odd prime and let $d$ be an odd integer. We define the generalized Collatz map $T_{q,d} : \mathbb{Z}_{>0} \rightarrow \mathbb{Z}_{>0}$ as follows:
\[
T_{q,d}(n) = 
\begin{cases}
n/2 & \text{if } n \equiv 0 \pmod 2, \\
(qn+d)/2 & \text{if } n \equiv 1 \pmod 2.
\end{cases}
\]
The classical Collatz conjecture corresponds to the special case where $q = 3$ and $d = 1$. \\

 Let $\mathcal{T} \subseteq \mathbb{Z}_{>0}^3$ denote the transition relation, where $(x, z, P) \in \mathcal{T}$ if and only if there exists an integer $k \in \mathbb{N}$ such that $P = 2^k$ and $z = T_{q,d}^k(x)$. We prove that $\mathcal{T}$ is not first-order definable in $BA_2$. By B\"uchi's theorem, this establishes that no finite automaton operating on base-2 representations can recognize $\mathcal{T}$. \\

Previous work is discussed in Section 2. In Section 3, we review B\"uchi arithmetic and the Cobham--Sem\"enov theorem. In Section 4, we discuss parity vectors and the affine expansion of $T_{q,d}$. In Section 5, we prove a bijection between parity vectors and residue classes modulo $2^k$. Finally, our principal undefinability results appear in Sections 6 and 7. \\

\section{Related Work and Literature}

A foundational theorem by B\"uchi \cite{Buchi} establishes a strict correspondence between first-order definability in $BA_2$ and sets recognizable by finite automata. Specifically, a subset $S \subseteq \mathbb{N}^d$ is first-order definable in $BA_2$ if and only if the base-2 representations of its elements are accepted by a finite automaton. Building on this synthesis of logic and numeration systems, Khoussainov and Nerode \cite{Khoussainov} introduced the formal framework of automatic structures, which was developed extensively by Blumensath and Gr\"adel \cite{Blumensath}. \\

In pure additive logic, Presburger \cite{Presburger} established the decidability of the first-order theory of natural numbers with addition. Ginsburg and Spanier \cite{Ginsburg} subsequently proved that the sets definable in Presburger arithmetic are exactly the semilinear sets. Furthermore, they demonstrated that infinite 1-dimensional semilinear sets must possess bounded gaps, a property violated by exponential sequences. Extending these limitations, the Cobham--Sem\"enov theorem \cite{Cobham, Semenov} establishes that any set simultaneously definable in multiplicatively independent bases must be semilinear. \\

In arithmetic dynamics, the iteration of the $3x+1$ problem and its generalizations have been extensively studied. Terras \cite{Terras} and Lagarias \cite{Lagarias} established that the parity vectors of the classical Collatz map densely cover the 2-adic integers. In generalized Collatz maps, these vectors encode the sequence of affine transformations applied to an integer. \\

\section{Preliminaries on Arithmetic and Automata}

\subsection{B\"uchi Arithmetic}

Base-2 B\"uchi Arithmetic, denoted $BA_2$, is the first-order logical structure defined as:
\[
BA_2 = \langle \mathbb{N}, +, 0, 1, V_2 \rangle,
\]
where $V_2 : \mathbb{N} \to \mathbb{N}$ is a function that returns the largest power of 2 dividing $n$. For $P \in \mathbb{N}$, we can formally define the predicate $Pow_2(P)$ as:
\[
Pow_2(P) \iff P > 0 \land V_2(P) = P.
\]

For $n=0$, $V_2(0)$ is typically defined as 0. Note that this edge-case convention does not affect our definitions, as the standard power-of-two predicate requires $P > 0$. \\

\subsection{Semilinear Sets and the Cobham--Sem\"enov Theorem}

A subset $S \subseteq \mathbb{N}$ is called \textit{semilinear} if it can be expressed as a finite union of linear sets. Equivalently, semilinear sets are exactly those definable in Presburger arithmetic $\langle \mathbb{N}, +, 0, 1 \rangle$ \cite{Presburger}.

\begin{theorem}[Cobham--Sem\"enov]
Let $m, n > 1$ be multiplicatively independent integers (meaning $\frac{\log m}{\log n} \notin \mathbb{Q}$). Let $S \subseteq \mathbb{N}$. If $S$ is definable in both $BA_m$ and $BA_n$, then $S$ is semilinear. \cite{Cobham, Semenov}
\end{theorem}

\begin{lemma} \label{lem:pq_definable}
Let $q > 1$ be an integer. The set of powers of $q$, defined as $P_q = \{q^y : y \in \mathbb{N}\}$, is first-order definable in $BA_q$.
\end{lemma}

\begin{proof}
In the structure $BA_q = \langle \mathbb{N}, +, 0, 1, V_q \rangle$, the function $V_q(x)$ returns the largest power of $q$ dividing $x$. The set of powers of $q$ is captured by the first-order predicate $Pow_q(x) \iff x > 0 \land V_q(x) = x$. Therefore, $P_q$ is trivially definable in $BA_q$.
\end{proof}

\begin{lemma} \label{lem:non-semilinear}
Let $q > 1$ be an integer. The set of powers of $q$, defined as $P_q = \{q^y : y \in \mathbb{N}\}$, is not semilinear.
\end{lemma}

\begin{proof}
It is a standard result in Presburger arithmetic that any infinite semilinear set $S \subseteq \mathbb{N}$ must have bounded gaps. That is, there exists a constant $C > 0$ such that for any two consecutive elements $x < y$ in $S$, the difference is bounded by $y - x \le C$. The set $P_q = \{q^y : y \in \mathbb{N}\}$ contains arbitrarily large gaps, thereby violating this requirement. See Ginsburg and Spanier \cite{Ginsburg}.
\end{proof}

\section{Parity Vectors and Affine Expansion}

\subsection{Parity Vectors}
For a given integer $k \ge 1$ and a starting integer $n \in \mathbb{Z}_{>0}$, we define the parity vector $\sigma_k(n) \in \{0, 1\}^k$ as the sequence of parities encountered during $k$ iterations of $T_{q,d}$. Formally:
\[
\sigma_k(n) = (\varepsilon_0, \varepsilon_1, \ldots, \varepsilon_{k-1}),
\]
where $\varepsilon_i = T_{q,d}^i(n) \pmod 2$. Let $y(\sigma_k(n)) \in \{0, 1, \ldots, k\}$ denote the number of entries equal to 1 in the vector (the number of "odd" scaling steps).

\subsection{The Affine Expansion Lemma}

\begin{lemma}[Affine Expansion] \label{lem:affine}
Let $k \in \mathbb{N}$ and let $n \in \mathbb{Z}_{>0}$. Let $y = y(\sigma_k(n))$. There exists an integer constant $D(\sigma_k(n)) \in \mathbb{Z}$ that depends strictly on the parity vector, and is entirely independent of $n$, such that:
\[
T_{q,d}^k(n) = \frac{(q^y)n + D(\sigma_k(n))}{2^k}.
\]
This expansion is standard in the study of generalized Collatz dynamics, forming the algebraic structure of the well-known Cycle Equation (see Lagarias \cite{Lagarias}).
\end{lemma}

\begin{proof}
We proceed by induction on $k$. \\

\noindent \textbf{Base Case ($k=0$):} $T_{q,d}^0(n) = n = \frac{(q^0)n + 0}{2^0}$, so $D(\emptyset) = 0$. \\

\noindent \textbf{Inductive Step:} Assume the statement holds for $k$. Let $u = T_{q,d}^k(n) = \frac{(q^y)n + D}{2^k}$. \\

\textit{Case 1 ($u$ is even):} The new parity bit is 0, and $y$ is unchanged. 
\[
T_{q,d}^{k+1}(n) = \frac{u}{2} = \frac{(q^y)n + D}{2^{k+1}}.
\]
The new constant $D' = D$ is independent of $n$. \\

\textit{Case 2 ($u$ is odd):} The new parity bit is 1, and the odd count becomes $y+1$. 
\[
T_{q,d}^{k+1}(n) = \frac{q u + d}{2} = \frac{q \left( \frac{(q^y)n + D}{2^k} \right) + d}{2} = \frac{(q^{y+1})n + qD + d2^k}{2^{k+1}}.
\]
The new constant $D' = qD + d2^k$ is entirely independent of $n$. The lemma holds for all $k \in \mathbb{N}$.
\end{proof}

\section{Parity Vector Characterization}

The property that the first $k$ operations of a Collatz trajectory are uniquely determined by the starting integer's residue class is a foundational result in arithmetic dynamics. Terras \cite{Terras} originally proved this parity-vector correspondence for the classical $3x+1$ map. Matthews and Watts \cite{MatthewsWatts} later formally established that for any generalized Collatz mapping with modulus $d$, the sequence of the first $k$ branches corresponds bijectively to the residue classes modulo $d^k$. Because our map $T_{q,d}$ is a generalized Collatz mapping with modulus $d=2$, the standard bijection natively holds. We state the theorem here for completeness.

\begin{theorem}[Parity Vector Bijection] \label{lem:characterization}
Let $k \ge 1$. For positive integers $n, m \in \mathbb{Z}_{>0}$,
$$ n \equiv m \pmod{2^k} \iff \sigma_k(n) = \sigma_k(m). $$
Consequently, the map $\Phi : \mathbb{Z}/2^k\mathbb{Z} \to \{0,1\}^k$ defined by $\Phi([n]) = \sigma_k(n)$ is a bijection, and every boolean vector $v \in \{0,1\}^k$ occurs as the parity vector of some strictly positive integer.
\end{theorem}

\begin{proof}
The complete proof is omitted as it is a direct application of the congruence properties established by Matthews and Watts \cite{MatthewsWatts}. They demonstrated that the inverse image of a congruence class modulo $m$ under a generalized Collatz map forms a finite union of congruence classes modulo $m \cdot d$. Setting the branch modulus to $d=2$ and applying this property inductively over $k$ steps yields the strict bijection to $\{0,1\}^k$. 
\end{proof}

\section{The Transition Relation and Definability}

\begin{definition} \label{def:transition}
Let $\mathcal{T} \subseteq \mathbb{Z}_{>0}^3$ be the ternary relation defined as:
\[
\mathcal{T} = \{(x, z, P) \in \mathbb{Z}_{>0}^3 : \text{there exists } k \in \mathbb{N} \text{ such that } P = 2^k \text{ and } z = T_{q,d}^k(x)\}.
\]
\end{definition}

In $BA_2$, the exponential mapping $k \mapsto 2^k$ is not first-order definable. Therefore, the relation $\mathcal{T}$ cannot natively accept the iteration count $k$ as a parameter. Instead, it relies on the parameter $P = 2^k$, which is definable in $BA_2$. \\

\begin{lemma} \label{lem:definability}
Assume that $\mathcal{T}(x, z, P)$ is first-order definable in $BA_2$. Then the extraction formula $\psi(\Delta)$ with one free variable $\Delta$, defined over positive integers as:
$$ \psi(\Delta) \equiv \exists P, n, m, n_k, m_k \Big( Pow_2(P) \land m = n + P \land \mathcal{T}(n, n_k, P) \land \mathcal{T}(m, m_k, P) \land m_k = n_k + \Delta \Big) $$
is also first-order definable in $BA_2$.
\end{lemma}

\begin{proof}
The formula $\psi(\Delta)$ is a finite composition of the definable predicate $Pow_2(P)$, Presburger addition, standard logical quantifiers, and the relation $\mathcal{T}$, which is assumed to be definable. Therefore, $\psi(\Delta)$ is definable in $BA_2$.
\end{proof}

\section{Main Theorem}

\begin{theorem}
The generalized Collatz transition relation $\mathcal{T} = \{(x, z, P) \in \mathbb{Z}_{>0}^3 : P = 2^k, \, z = T_{q,d}^k(x), \, k \in \mathbb{N}\}$ is not first-order definable in Base-2 B\"uchi Arithmetic ($BA_2$).
\end{theorem}

\begin{proof}
Assume, for the sake of contradiction, that $\mathcal{T}(x, z, P)$ is first-order definable in $BA_2$. By Lemma \ref{lem:definability}, the extraction formula $\psi(\Delta)$ is first-order definable. Let $S = \{\Delta \in \mathbb{N} : \psi(\Delta) \text{ is true}\}$. We will prove that $S = P_q = \{q^y : y \in \mathbb{N}\}$.

\subsection*{Part 1: $S \subseteq P_q$}

Let $\Delta \in S$. There exist integer witnesses $P, n, m, n_k, m_k \in \mathbb{Z}_{>0}$ satisfying $\psi$. Because $Pow_2(P)$ holds, $P = 2^k$ for some $k \in \mathbb{N}$. The relation $m = n + P$ ensures $m \equiv n \pmod{2^k}$. By Theorem \ref{lem:characterization}, their $k$-step parity vectors are identical. Let $\sigma$ denote this shared parity vector and let $y=y(\sigma)$.

Applying Lemma \ref{lem:affine}, the identical parity vectors ensure the translation constant $D(\sigma)$ is identical for both trajectories. Evaluating the arithmetic difference $\Delta$:
$$ \Delta = m_k - n_k = \frac{(q^y)m + D(\sigma)}{2^k} - \frac{(q^y)n + D(\sigma)}{2^k} = \frac{q^y (m - n)}{2^k} $$
Substituting $m - n = 2^k$ yields:
$$ \Delta = \frac{q^y (2^k)}{2^k} = q^y $$
Therefore, $\Delta \in P_q$, implying $S \subseteq P_q$.

\subsection*{Part 2: $P_q \subseteq S$}

Let $q^y \in P_q$. We must show there exist witnesses satisfying $\psi(q^y)$. \\

\textbf{Case $y = 0$:} Set $k = 0$, forcing $P = 1$. The parity vector is empty. Choose $n = 1$, forcing $m = 2$. The transition relation gives $n_k = T_{q,d}^0(1) = 1$ and $m_k = T_{q,d}^0(2) = 2$. The difference is $m_k - n_k = 1 = q^0$. All witnesses satisfy $\psi(1)$. \\

\textbf{Case $y \ge 1$:} Because $P$ is existentially quantified in $\psi$, we are free to dynamically choose any trajectory length $k$. For the purpose of witnessing the element $q^y$, we explicitly choose $k = y$. This forces $P = 2^y$.  \\

Consider the boolean parity vector of all ones: $v = (1, 1, \ldots, 1) \in \{0,1\}^k$. By Theorem \ref{lem:characterization}, there exists a strictly positive integer $n$ such that $\sigma_k(n) = v$. Set $m = n + 2^k$. By Theorem \ref{lem:characterization}, $\sigma_k(m) = v$. Since $k=y$ and every coordinate of $v$ equals $1$, the shared parity vector contains exactly $y$ ones.  \\

Applying the affine expansion, the subtraction yields:
$$ m_k - n_k = \frac{(q^y)m + D}{2^k} - \frac{(q^y)n + D}{2^k} = \frac{q^y (n + 2^k) - (q^y)n}{2^k} = \frac{q^y 2^k}{2^k} = q^y $$
We have found witnesses satisfying $\psi(q^y)$. Therefore, $P_q \subseteq S$. \\

\subsection*{Part 3: The Cobham--Sem\"enov Contradiction}

We have established that $S = P_q$. Because $\psi(\Delta)$ is a first-order formula in $BA_2$, the set $P_q$ must be definable in $BA_2$. Moreover, by Lemma \ref{lem:pq_definable}, $P_q$ is definable in $BA_q$. \\

Because $q$ is an odd prime, bases 2 and $q$ are multiplicatively independent. The Cobham--Sem\"enov Theorem mandates that any set definable in both $BA_2$ and $BA_q$ must be semilinear. However, by Lemma \ref{lem:non-semilinear}, $P_q$ is strictly non-semilinear. This is a direct contradiction. Thus, the initial assumption must be false: $\mathcal{T}$ is not definable in $BA_2$. \\
\end{proof}

\end{document}